\documentclass[a4paper]{article}
\usepackage[affil-it]{authblk}

\usepackage{xcolor,amsthm,amsmath,amssymb,xspace}
\usepackage{epstopdf}% To incorporate .eps illustrations using PDFLaTeX, etc.
\usepackage[caption=false]{subfig}% Support for small, `sub' figures and tables

\usepackage[numbers,sort&compress]{natbib}% Citation support using natbib.sty
\bibpunct[, ]{[}{]}{,}{n}{,}{,}% Citation support using natbib.sty
\renewcommand\bibfont{\fontsize{10}{12}\selectfont}% Bibliography support using natbib.sty
\makeatletter% @ becomes a letter
\def\NAT@def@citea{\def\@citea{\NAT@separator}}% Suppress spaces between citations using natbib.sty
\makeatother% @ becomes a symbol again

\theoremstyle{plain}% Theorem-like structures provided by amsthm.sty
\newtheorem{theorem}{Theorem}[section]

\theoremstyle{definition}

\theoremstyle{remark}

\DeclareMathOperator{\R}{\mathbb{R}}

\usepackage{tikz}
\usetikzlibrary{matrix,chains,positioning,decorations.pathreplacing,arrows}

\begin{document}

%\articletype{ARTICLE TEMPLATE}% Specify the article type or omit as appropriate

\title{Best  free knot linear spline approximation and its  application to neural networks}

\author{
{D.~K. Pham and V. Peiris and N. Sukhorukova\thanks{nsukhorukova@swin.edu.au}}
%\affiliation{\textsuperscript{a}Swinburne University of Technology, Australia; \textsuperscript{b}Curtin University, Australia}
}

\maketitle

\begin{abstract}
The problem of fixed knot approximation is convex and there are several efficient approaches to solve this problem,  yet, when the knots joining the affine parts are also variable, finding conditions for a best Chebyshev approximation remains an open problem. It was noticed before that piecewise linear approximation with free knots is equivalent to neural network approximation with piecewise linear activation functions (for example ReLU).  In this paper, we demonstrate that in the case of one internal free knot, the problem of linear spline approximation can be reformulated as a mixed-integer linear programming problem and solved efficiently using, for example, a branch and bound type method. We also present a new sufficient optimality condition for a one free knot piecewise linear approximation.  The results of numerical experiments are provided.
\end{abstract}


\section{Introduction}
The problem of approximating a continuous function by a continuous piecewise polynomial  functions (splines) has been studied for several  decades~\cite{Nurnberger, Schum}. When the points of joining the polynomial pieces (knots) are also variables, finding conditions for a best Chebyshev approximation remains an open
problem~\cite{Nurnberger}. This problem was identified as one the of the most important problems of modern approximation~\cite{OpenProblems}. There are a number of theoretical results on free knot polynomial spline approximation: necessary optimality conditions~\cite{Nurnberger, NurnSchuSomm, SukhUgon2017Transactions, CrouSukhUgon}, sufficient optimality conditions~\cite{Nurnberger, NurnSchuSomm} (just to name a few). As in most non-convex problems, necessary and sufficient optimality conditions do not coincide and it is not known if the conditions can be improved.

In this paper, we consider a special class of free knot polynomial spline approximation, where the degree of the polynomial pieces does not exceed one (that is, approximation by free knot linear splines).  This problem is simpler than the  general free knot polynomial spline approximation, but it is also equivalent to approximation by a neural network with one hidden layer and ReLU (Leaky ReLU) activation function~\cite{Chen}, where the number of the nodes in the hidden layer does not exceed the number of distinct knots. All the results are for Chebyshev (uniform) norm. In this paper we are focusing on linear splines with only one knot. At first glance, this problem is simple, but this is  the first step to create efficient Remez-type algorithms for free knot spline approximation. Moreover, this problem has many practical applications, including data science and deep learning. 

Artificial neural networks are very powerful and popular approximation techniques with many practical applications, including image and sound recognition, financial fraud detection, partial differential equation, fluid dynamics and many others, see~\cite{Goodfellow2016} and references therein. The list of applications where deep learning demonstrated its efficiency is growing.  Neural networks can handle models with a large number of variables: complex models are decomposed into several simple optimisation problems that can be solved efficiently.  

Deep learning is based on artificial neural networks where the network has a specific structure: the input layer, one or more hidden layers and the output layer. In this paper, we consider artificial neural networks with a single hidden layer and therefore, in the context of this paper, deep learning algorithms and artificial neural networks algorithms refer to the same techniques.

The goal of this paper is to start with linear spline approximation with only one free internal knot.  The contribution of this paper is as follows.
\begin{enumerate}
    \item We reformulate one free knot linear spline problems as mixed-integer linear programming problems (MILP), whose solutions can be obtained using a branch and bound type of method.
    \item We develop strong sufficient optimality conditions for one free knot linear spline approximation problems that are stronger than the existing general case conditions~\cite{Nurnberger}.
\end{enumerate}
We also present the results of numerical experiments and compare them with the results obtained by neural network approximations.
There have been a number of attempts to explain the efficiency of neural networks, based on their optimisation properties~\cite{Goodfellow2016, Sun2020OptimDeepLearning, Vidal2017},  but there are still many open problems remain.

The paper is organised as follows. In Section~\ref{sec:Prelim}, we provide some preliminaries on free knot linear approximation, neural network approximation and their connections. Section~\ref{sec:MILPP} discusses the reformulation of the one free knot linear spline problem as a mixed-integer linear programming problem. Sufficient optimality conditions are presented in Section~\ref{sec:SufficientCond} and in Section~\ref{sec:Experiments}, we discuss the results of our numerical experiments. Finally, we summarise our work and comment on future research directions in Section~\ref{sec:Conclusion}.

\section{Preliminaries} \label{sec:Prelim}
\subsection{Free knot linear approximation in Chebyshev norm}
 Assume that a continuous univariate function $f(t)$ is to be approximated by a continuous piecewise linear function (linear spline) in a bounded interval $[c,d]$. The points $c$ and $d$ are called the {\it external knots}. The external knots are fixed. The points where the piecewise linear function switches from one linear piece to the next one are called the {\it internal knots}. If $n$ is the number of subintervals within $[c,d]$, then the number of internal knots is $n-1$:
 \begin{equation*}   c=\theta_0\leq\theta_1\leq\theta_2\leq\dots\leq\theta_{n-1}\leq\theta_n=d.
 \end{equation*}
 
A free knot linear spline approximation problem (that is, the internal knots are free and part of the decision variables) can be formulated as follows:
\begin{equation}\label{eq:free_cont}
 {\text{minimise}}\sup_{t\in [c,d]} \left| a_0+\sum_{i=1}^{n}a_i\max\{0,t-\theta_{i-1}\}- f(t)\right|, ~{\text{subject to}}~ X\in \R^{2n}
\end{equation}
where $n$ is the number of subintervals, $\theta_0=c$, $\theta_n=d$, $A=(a_0, a_1,\dots,a_n)$ are called the spline coefficients (spline parameters),  $X=(a_0,a_1,\dots,a_n,\theta_1,\dots,\theta_{n-1})$ are the decision variables.

In most practical problems, one needs to discretise the the domain and the problem is as follows:
\begin{equation}\label{eq:free_disc}
 {\text{minimise}}\max_{t_j,  j=1,\dots,N} \left| a_0+\sum_{i=1}^{n}a_i\max\{0,t_j-\theta_{i-1}\}- f(t)\right|, ~{\text{subject to}}~ X\in \R^{2n}
\end{equation}
where $N$ is the number of discretisation points, $t_j\in[c,d]$ for $j=1,\dots,N$,  $n$ is the number of subintervals, $\theta_0=c$, $\theta_n=d$, $X=(a_0,a_1,\dots,a_n,\theta_1,\dots,\theta_{n-1})$ are the decision variables. Most numerical algorithms are based on discretisation and therefore this is the main target of our paper.

If the internal knots are fixed (fixed knot spline approximation) the corresponding optimisation problem is convex and there are a number of efficient methods to solve this problem. If, however, the knots are part of the decision variables, the problem is still open. 

In this paper, we consider free knot linear splines with only one internal knot, that is, $n=2$. These problems are essential for further extension to a larger number of subintervals. Our approach is based on a reformulation of the original problem as a pair of mixed-integer linear programming problems, whose solution can be obtained by a branch and bound type of optimisation methods, available in most modern optimisaton packages (section~\ref{sec:MILPP}).  In this paper, we also develop sufficient optimality conditions for linear spline approximation with one free internal knot (section~\ref{sec:SufficientCond}). These conditions are used to verify optimality.

Apart from being a very important open problem, there are a number of practical applications for free knot linear spline approximation, including data approximation, deep learning and many others. In the next subsection, we provide a short introduction to neural networks and demonstrate the connections between free knot linear spline approximation and neural network approximation. In particular, these approaches are equivalent in the case of univariate function approximation. 

Uniform (Chebyshev) approximation problems and their optimality conditions (necessary, sufficient) are based on the number of alternating points, that is points whose absolute deviation is maximal and the signs of deviation are alternating~\cite{Nurnberger,Schum,NurnSchuSomm,SukhUgon2017Transactions,CrouSukhUgon}, just to name a few. The notion of alternating points and alternating sequences is central for uniform approximation and it is not specific for polynomial splines (free or fixed knots).  

The origins of uniform approximation start with P.~Chebyshev~\cite{Cheb}.
\begin{theorem}\label{thm:cheb}(Chebyshev) A necessary and sufficient optimality condition for a polynomial of degree~$n$ is the existence of $n+2$ alternating points. 
\end{theorem}
In particular, if the degree of the polynomial is one (linear approximation), the necessary and sufficient optimality condition is the existence of three alternating points.

\subsection{Neural networks for approximations}

Neural networks is a popular tool in the modern area of Machine learning and Artificial Intelligence. Deep learning is just a subclass of artificial neural networks that has many practical applications, including data analysis, signal and image processing and many others~\cite{Goodfellow2016,Sun2020OptimDeepLearning}.  This popularity is due to the efficiency of approximation models (and methods) behind deep learning. Moreover, deep learning can be used as a purely approximation tool for univariate and multivariate functions.  The goal of deep learning is to optimise weights in the network and therefore, this problem can be seen as an approximation problem treated using modern optimisation tools. In the corresponding optimisation models, the objective function represents the deviation (error) of the approximation from the original function (which may be a continuous function or a function whose values are known only at some discretisation points). The deviation function (called ``loss function'') can be the sum of the squares of the deviations (least squares based models), the maximum of the absolute values (uniform or Chebyshev based models), the sum of the absolute values (Manhattan-distance based models), etc.

Deep learning is based on solid mathematical modelling established in~\cite{Cybenko,Hornik1991,Pinkus1993,Pinkus1999},  but the origins come to the work of A.~Kolmogorov and his student V.~Arnold. The  celebrated Kolmogorov-Arnold Theorem~\cite{Kolmogorov1957, Arnold1963} is an attempt to solve the 13th problem of Hilbert.

It is customary to choose the mean least squares loss functions in neural networks. There are several reasons for this. First of all, this model involves minimising a smooth quadratic objective function, and therefore basic optimisation techniques such as the gradient descent can be used. Second, the least squares loss function also fares well with the assumption that the errors outside the discretisation points (where the value of the function is known) have a normal distribution. At the same time, the approximation results established in~\cite{Cybenko,Hornik1991,Pinkus1993,Pinkus1999} rely on the Chebyshev (uniform) convergence, which is a stronger result. Therefore, there are a number of interesting results where the loss function is Chebyshev norm based~\cite{PeirisRoshchinaSukhorukova2023}.

% \begin{figure}
% \centering
% %\includegraphics[width=0.45\textwidth]{neural-network-no-hidden-upd.pdf}
% \includegraphics[width=0.45\textwidth]{MultiLayerNeuralNetwork.pdf}
% \label{fig:1}
% \caption{Neural Network: one hidden layer}
% \end{figure} 

\subsection{ReLU networks and linear splines}
\begin{figure}
    \centering
    \includegraphics[width=0.5\textwidth]{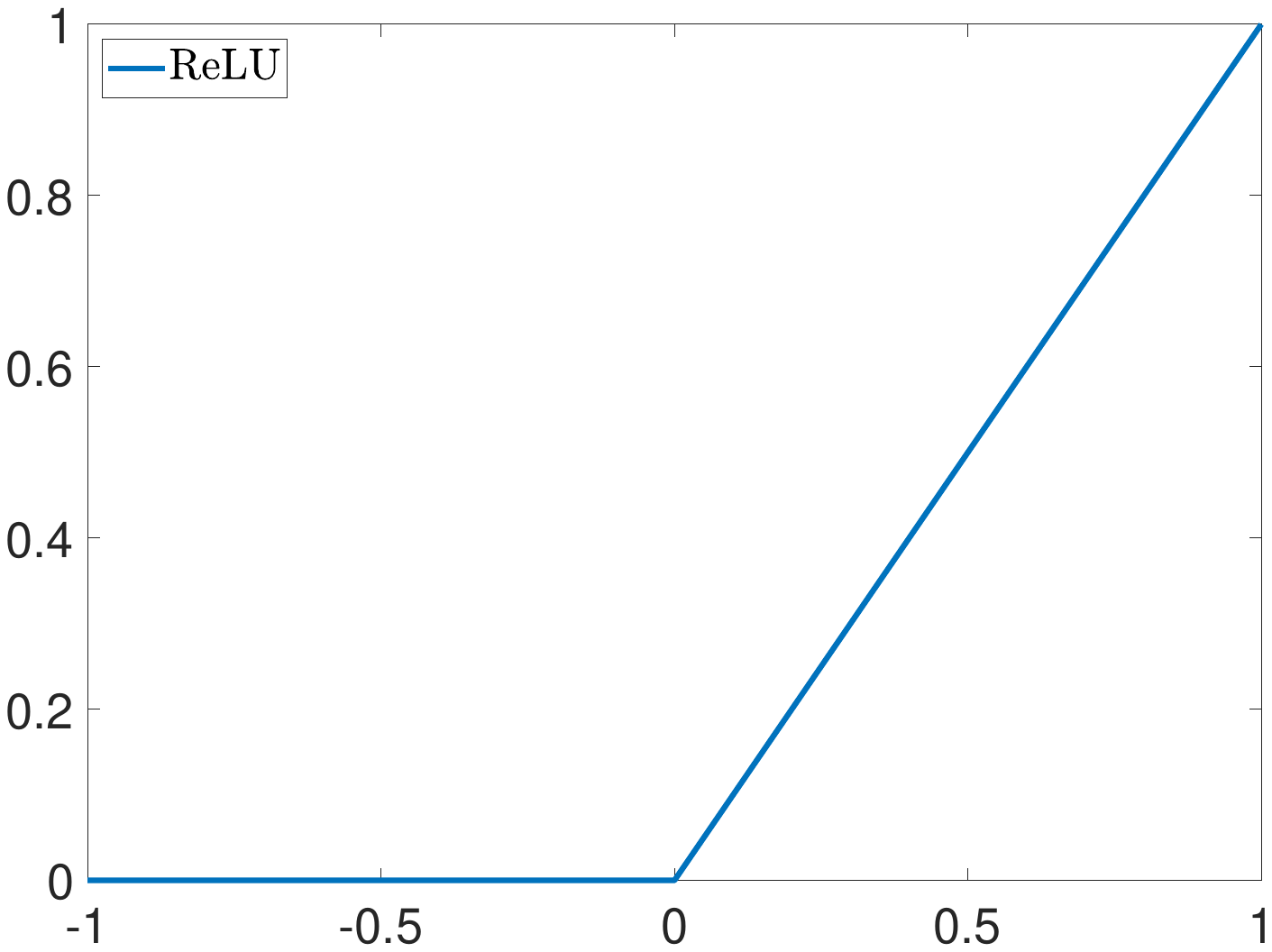}
    \caption{ReLU function in $[-1,1]$}
    \label{fig:ReLU}
\end{figure}

The ReLU (Rectified linear unit) function is defined as follows:
$$g(a) = \max\{0,a\} = \begin{cases}
		  a, & \text{if} \quad a > 0,\\
            0, & \text{otherwise,}
		 \end{cases}$$
where $a$ is the input from a given interval $[a_{min}, a_{max}]$. Figure~\ref{fig:ReLU} shows ReLU function defined in $[-1,1].$

The function $g(a)$ has a discontinuity at $0$ and therefore, it can be considered as a linear spline with one internal knot at $0$.
%The end points are $g(a_{min}) = 0$ and $g(a_{max}) = a$. 

Let $z$ be a single scalar input for the same function $g$ above and $z=wx+b$
%inner node in the hidden layer of the neural network defined as 
where $x$ is the input from the interval $[x_{min}, x_{max}]$, $w$ and $b$ are decision variables. First assume that $w\ne 0$, then, if $w>0$, we have:
\begin{align*} 
g(z) = g(wx+b)  &=  \max\{0,wx+b\} \\ 
                &=  \begin{cases}
                        wx+b, & \text{if} \quad wx+b > 0,\\
                        0, & \text{otherwise.}
		          \end{cases}\\
                &=  \begin{cases}
                        wx+b, & \text{if} \quad x > \frac{-b}{w},\\
                        0, & \text{otherwise.}
		          \end{cases}
\end{align*}

In a similar manner, for $w<0$, we have
\begin{align*} 
g(z) = g(wx+b)  &=  \max\{0,wx+b\} \\                 		         
                &=  \begin{cases}
                        wx+b, & \text{if} \quad x < \frac{-b}{w},\\
                        0, & \text{otherwise.}
		          \end{cases}
\end{align*}

This means that there is a discontinuity at $\frac{-b}{w}$ and hence, there is an internal knot located at $x=\frac{-b}{w}$. 
%The location of the knot depends on the sign of the decision variables. In particular, if $w$ and $b$ both have the same sign (negative or positive), then the location of the knot is $\frac{-b}{w}$. If they have different signs, then the location of the knot is $\frac{b}{w}$. 
A single linear piece will result in both cases if the knot is located outside of the considered domain.
%The end points will be shifted to $[x_{min}-\frac{-b}{w}, x_{max}-\frac{-b}{w}]$. 
Finally, assume that $w = 0,$ then 
\begin{align*} 
g(z) = g(b)  &=  \max\{0,b\} \\ 
                &=  \begin{cases}
                        b, & \text{if} \quad b > 0,\\
                        0, & \text{otherwise.}
		          \end{cases}
\end{align*}
where the output depends on the value of $b$ regardless of the input $x$. In this case, $g(z)$ consists of only one linear piece. %parallel to the x-axis. 

To summarise, each node in the hidden layer neural network approximation with ReLU activation function leads to a knot in the equivalent linear spline approximation.  In the case when $w\ne 0$, the location of the knot is $\theta=-\frac{b}{w}$. In the case when $w=0$ the corresponding component is a single linear piece, which is equivalent to allocating a knot outside the approximation interval~$[c,d]$. 

A natural step now is to extend the result to the neural networks with more than one node in their hidden layers. Intuitively, we expect that each node leads to a knot, but some knots may coincide and some may ``move'' outside the approximation interval~$[c,d]$ and therefore the total number of internal knots in the linear spline approximation does not exceed the total number of nodes in the hidden layer. In the rest of this section, we present a formal explanation of this observation.  

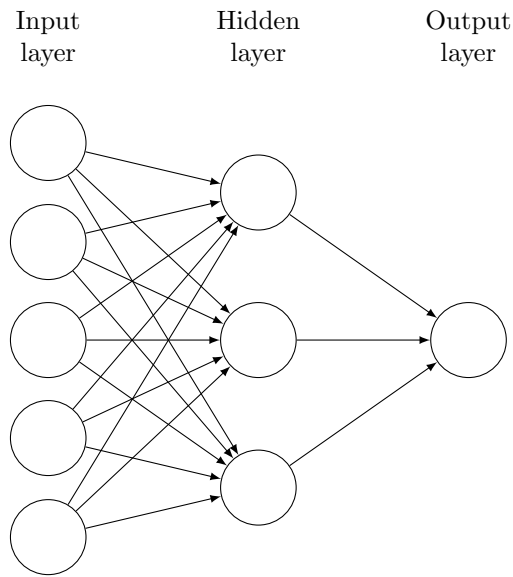
\begin{figure}
    \centering
    \begin{tikzpicture}[
             % define styles 
             clear/.style={ 
                 draw=none,
                 fill=none
             },
             net/.style={
                 matrix of nodes,
                 nodes={ draw, circle, inner sep=10pt },
                 nodes in empty cells,
                 column sep=0.3cm,
                 row sep=-10pt
             },
             >=latex
        ][scale = 0.1pt]
        % define matrix mat to hold nodes
        % using net as default style for cells
        \matrix[net] (mat)
        {
        % Define layer headings
        |[clear]| \parbox{1.3cm}{\centering Input\\layer} 
            & |[clear]| \parbox{1.3cm}{\centering Hidden\\layer} 
            & |[clear]| \parbox{1.3cm}{\centering Output\\layer} \\
                 
          & |[clear]|        & |[clear]| \\
        |[clear]|         &  & |[clear]| \\
          & |[clear]|        & |[clear]| \\
        |[clear]|         & |[clear]|        & |[clear]|  \\
          &  &  \\
        |[clear]|         & |[clear]|        & |[clear]|   \\
          & |[clear]|        & |[clear]| \\
        |[clear]|         &  & |[clear]| \\
          & |[clear]|        & |[clear]| \\ 
        };
        % lines from a_{i}^{0} to each a_{j}^{1}
        \foreach \ai in {2,4,...,10} {
            \foreach \aii in {3,6,9}
                \draw[->] (mat-\ai-1) -- (mat-\aii-2);
                }
        % lines from a_{i}^{1} to a_{0}^{2}
        \foreach \ai in {3,6,9}
          \draw[->] (mat-\ai-2) -- (mat-6-3);
    \end{tikzpicture}
    \caption{Neural Network: one hidden layer}
    \label{fig:NN}
\end{figure}

Consider a neural network with one hidden layer where the activation function defined on the hidden nodes is ReLU and the output layer consists of just one node. This is a multilayered network (Figure~\ref{fig:NN}) whose final output is a composition of linear functions with the activation function ReLU. 

Let $w^{(i)}_j,$ $i=1,2, j=1,\ldots,n$ be the weights of the network, $b^{(i)}_j,$ $i=1, j=1,\ldots,n$ are the bias terms assigned to the nodes in the hidden layer, $b^{(2)}$ is the bias term assigned to the output of the network.

If the input for the inner nodes in the hidden layer are denoted as $z_1, z_2, \ldots, z_n$ where $n$ is the number of inner nodes in the hidden layer, then the output of hidden nodes (input for the output layer) can be defined as follows:
\begin{align*}
    g(z_1) &= g(w^{(1)}_1 x + b^{(1)}_1),\\
    g(z_2) &= g(w^{(1)}_2 x + b^{(1)}_2),\\
    &\vdots\\
    g(z_n) &= g(w^{(1)}_n x + b^{(1)}_n).
\end{align*}

If the neural network output is denoted by $y$, then $y$ is an affine transformation of inner nodes:
\begin{align*}
    y(x) &= w^{(2)}_1 g(z_1) + w^{(2)}_2 g(z_2) + \ldots + w^{(2)}_n g(z_n) + b^{(2)},\\
    y(x) &= w^{(2)}_1 g(w^{(1)}_1 x + b^{(1)}_1) + w^{(2)}_2 g(w^{(1)}_2 x + b^{(1)}_2) + \ldots + w^{(2)}_n g(w^{(1)}_n x + b^{(1)}_n) + b^{(2)}.
\end{align*}
The union of all the knots created by the inner nodes represents the set of knots of the function $y$. As $y$ consists of a linear combination of $z_i$'s, one can see that $y$ is a linear spline.

A neural network (one hidden layer) with ReLU activation is a linear spline and the number of the nodes in the hidden layer does not exceed the number of distinct knots of the linear spline~\cite{Chen, hansson2017}. In this case, finding the internal knots is easier due to the natural architecture of the network itself. The weights and biases of the network are adjustable parameters decided by the training (optimisation) of the network. Hence, the location of the knots of the approximation computed by ReLU network with one hidden layer is decided in the optimisation of the weight space of the network.

%    \begin{figure}\label{fig:2}
 %   \centering
%\includegraphics[width=0.5\textwidth]{MultiLayerNeuralNetworkVector.pdf}
%    \caption{Neural Network (Wikipedia picture: creative commons)}
%\end{figure} 

\section{Mixed integer linear programming approach}\label{sec:MILPP}

In the case of two subintervals (one free knot), the solution  can be found by solving two mixed-integer linear programming problems. In the first problem, the optimal linear spline is the maximum of two affine pieces, while in the second problem the optimal linear spline is the minimum of two affine pieces. These pieces may be represented by the same affine function (single line).  

The main advantage of this approach is that the knot does not need to be calculated explicitly.  Indeed, since there are at most two subintervals, the corresponding piecewise linear approximation is either a maximum or a minimum of two affine functions. For each case the problem can be formulated a mixed-integer linear programming problem. After solving both problems,   the one with minimal optimal value is the solution to the original free knot linear spline problem.

In the rest of this section, we work with the discretised problem~(\ref{eq:free_disc}). Assume that $a_1$ and $b_1$ are the slope and the intercept of the affine piece in the first subinterval, respectively and $a_2$ and $b_2$ are the slope and the intercept of the affine piece in the second subinterval, respectively.

\subsection{Maximum problem}

First, for each discretisation point $t_i,$ $i=1,\dots,N$ we introduce a new variable 
\begin{equation}\label{eq:maxci} 
c_i=\max_{i}\{a_1 t_i+b_1,a_2 t_i+b_2\},~i=1,\dots,N,\end{equation} 
where $N$ is the number of discretisation points and may be large. The objective is to minimise the absolute deviation $z$, subject to the following constraints:
\begin{align}\label{eq:constraint_dev1}
f(t_i)-c_i\leq z,~i=1,\dots,N,\\
\label{eq:constraint_dev2}
c_i-f(t_i)\leq z,~i=1,\dots,N.
\end{align}

Due to~(\ref{eq:maxci}), we have the following equations:
\begin{align}\label{eq:max_gr1_1}
a_1t_i+b_1\leq c_i,~i=1,\dots,N,  \\
\label{eq:max_gr1_2}
a_2t_i+b_2\leq c_i,~i=1,\dots N
\end{align}
and for every $i$, at least one of the inequalities  has to be satisfied as equality. This is where we have to introduce a binary variable. This can be achieved by requiring that for each group~$i$, at least one of the following reverse inequalities is satisfied:
\begin{align*}
a_1t_i+b_1\geq c_i,~i=1,\dots,N,    \\
a_2t_i+b_2\geq c_i,~i=1,\dots N.
\end{align*}
For each group, introduce a binary variable $z_i$. Also consider a larger positive parameter $M$ (big-$M$, fixed value). Then the requirement for at least one of the inequalities in each group holds can be expressed as follows:
\begin{align}
c_i-(a_1t_i+b_1)\leq Mz_i,~i=1,\dots,N,    \\
c_i-(a_2t_i+b_2)\leq M(1-z_i),~i=1,\dots N,\\
z_i\in\{0,1\},~i=1,\dots,N.\label{eq:max_gr3}
\end{align}
Finally, the goal is to minimise~$z$ subject to~(\ref{eq:constraint_dev1})-(\ref{eq:max_gr3}).

\subsection{Minimum problem}

Similar to the maximum problem, introduce a new variable 
\begin{equation}\label{eq:mindi} 
d_i=\min_{i}\{a_1 t_i+b_1,a_2 t_i+b_2\},~i=1,\dots,N,\end{equation} 
where $N$ is the number of discretisation points. The objective is to minimise the absolute deviation $y$, subject to 
\begin{align}\label{eq:constraint_dev3}
f(t_i)-d_i\leq y,~i=1,\dots,N,\\
\label{eq:constraint_dev4}
d_i-f(t_i)\leq y,~i=1,\dots,N.
\end{align}

Due to~(\ref{eq:mindi}), we have the following equations:
\begin{align}\label{eq:min_gr1_1}
a_1t_i+b_1\geq d_i,~i=1,\dots,N,  \\
\label{eq:min_gr1_2}
a_2t_i+b_2\geq d_i,~i=1,\dots N
\end{align}
and for every $i$, at least one of the inequalities  has to be satisfied as equality. For each group~$i$, at least one of the following reverse inequalities is satisfied:
\begin{align*}
a_1t_i+b_1-d_i\leq 0,~i=1,\dots,N,    \\
a_2t_i+b_2-d_i\leq 0,~i=1,\dots N.
\end{align*}
For each group, introduce a binary variable $y_i$ and a large positive number~$M$. Then the final block of inequalities is as follows:
\begin{align}
a_1t_i+b_1-d_i\leq M(1-y_i),~i=1,\dots,N,    \\
a_2t_i+b_2-d_i\leq My_i,~i=1,\dots N,\\
y_i\in\{0,1\}\label{eq:min_gr3}
\end{align}
Finally, the problem is to minimise~$y$ subject to~(\ref{eq:constraint_dev3})-(\ref{eq:min_gr3}).

\section{Sufficient optimality conditions}\label{sec:SufficientCond}

In this section, we develop sufficient optimality conditions for one free knot linear spline approximation. We have to treat two cases separately: 
\begin{enumerate}
    \item the optimal spline contains two distinct linear segments;
    \item the optimal spline contains a single linear segment.
\end{enumerate}

\begin{theorem}\label{thm:sufficient}
    The following two conditions are sufficient optimality conditions for free knot linear spline approximation in Chebyshev norm with one internal knot.
    \begin{enumerate}
        \item The linear spline consists of two distinct linear pieces has at least three alternating points in each subinterval.
        \item The linear spline consists of a single linear piece and has at least four alternating points in the whole interval.
    \end{enumerate}
\end{theorem}
Proof: 
Let $f(t)$ be a continuous function which needs to be approximated in the interval~$[c,d]$.
\begin{enumerate}
    \item Assume that there are two distinct linear pieces and at least three alternating points in each subinterval. Due to theorem~\ref{thm:cheb}, at each subinterval the approximation is optimal and the maximal deviation can not be improved by changing the coefficients of the spline. The change of the location of the knot also can not improve the maximal deviation, since by increasing the length of one of the subintervals the maximal absolute deviation is non-decreasing and therefore the linear spline is optimal.
    \item Now assume that the spline consists of only one linear piece $l(t)$ and there are four alternating points: $t_1, t_2, t_3$ and $t_4$. Due to theorem~\ref{thm:cheb}, the current spline is optimal, unless the optimal spline consists of two linear pieces.

    Assume that the maximal deviation can be improved  by splitting the original interval~$[c,d]$ into two subintervals: $[c,\theta]$ and $[\theta,d]$. The corresponding linear pieces are $l_1(t)$ and $l_2(t)$:
    $$s(t)=
    \begin{cases}
        l_1(t),& t\in[c,\theta]\\
        l_2(t),& t\in(\theta,d].
    \end{cases}
        $$
    It is clear that the knot $\theta\in(t_2,t_3)$, since otherwise one of the intervals contains three alternating points and therefore the maximal deviation can not be improved. 
    
    Without loss of generality, assume that $$0<\delta=f(t_1)-l(t_1).$$ Since $s(t)$ provides a smaller absolute deviation, 
    $$0<\delta=f(t_1)-l(t_1)>f(t_1)-s(t_1)=f(t_1)-l_1(t_1)$$
    and $$0>-\delta=f(t_2)-l(t_2)<f(t_2)-s(t_2)=f(t_2)-l_1(t_2).$$
    Therefore, $s(t_1)=l_1(t_1)>l(t_1)$ and $s(t_2)=l_1(t_2)<l(t_2)$. Hence, 
    \begin{equation}\label{eq:theta1}
        s(\theta)=l_1(\theta)<l(\theta). 
    \end{equation}
    In a similar manner, the analysis of the behaviour of the linear pieces $l_2(t)$ and $l(t)$ yields
    \begin{equation}\label{eq:theta2}
        s(\theta)=l_2(\theta)>l(\theta). 
    \end{equation}
    Combining~(\ref{eq:theta1}) and (\ref{eq:theta2}), conclude that $s(t)$ is discontinuous at $\theta$, which contradicts to the requirement for the optimal linear spline to be continuous at its knot.      
\end{enumerate}

\hskip300pt $\square$
\section{Numerical experiments} \label{sec:Experiments}
All the numerical experiments are performed on the interval $[-1,1]$ and the discretisation step is~$h=10^{-3}$. 
We test our method on five different functions:
\begin{enumerate}
\item $f_1(t)=\sqrt{|t|}$; this function is nonsmooth and non-Lipschitz.
\item $f_2(t)=\sqrt{|t-0.75|}$; this function is similar to~$f_1(t)$, but it is non-symmetric.
\item $f_3(t)=\sin(2\pi t)$; this function is periodic and oscillating.
\item $f_4(t)=t^3-3t^2+2$; this is a cubic function.  The experiments with this function are interesting, since this is an example, where neural network was especially inaccurate, despite the fact that this is a smooth function without any abrupt changes. 
\item $f_5=1/(t^{25}+0.5)$; this is a very complex function for approximating by a continuous piecewise linear function with only two linear pieces. The structure of the approximation drastically changes when the the discretisation step is changing.
\end{enumerate}

We approximate above functions by one free knot piecewise linear functions and by one hidden layer neural network with ReLU activation function. The hidden layer consists of just one node which represents the internal knot of the piecewise linear approximation. The loss function is based on the uniform (Chebyshev) approximation and the optimiser, ADAM and ADAMAX (a certain adaptation of ADAM to uniform approximation) used to find the parameters of the network. 

\begin{figure}
    \centering
    \includegraphics[width = 0.3\textwidth]{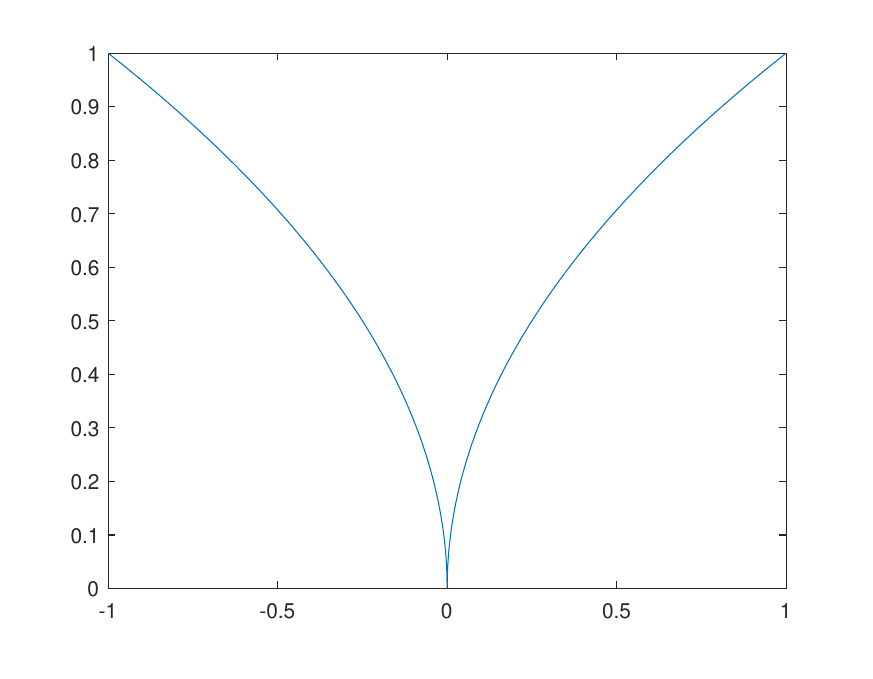}
    \includegraphics[width = 0.33\textwidth]{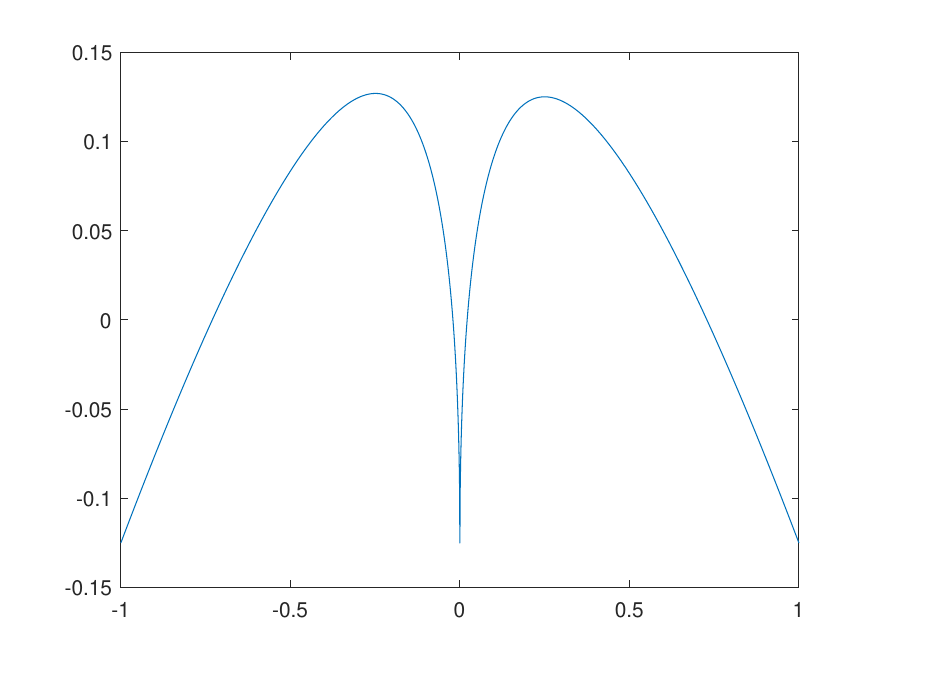}
    \includegraphics[width = 0.29\textwidth]{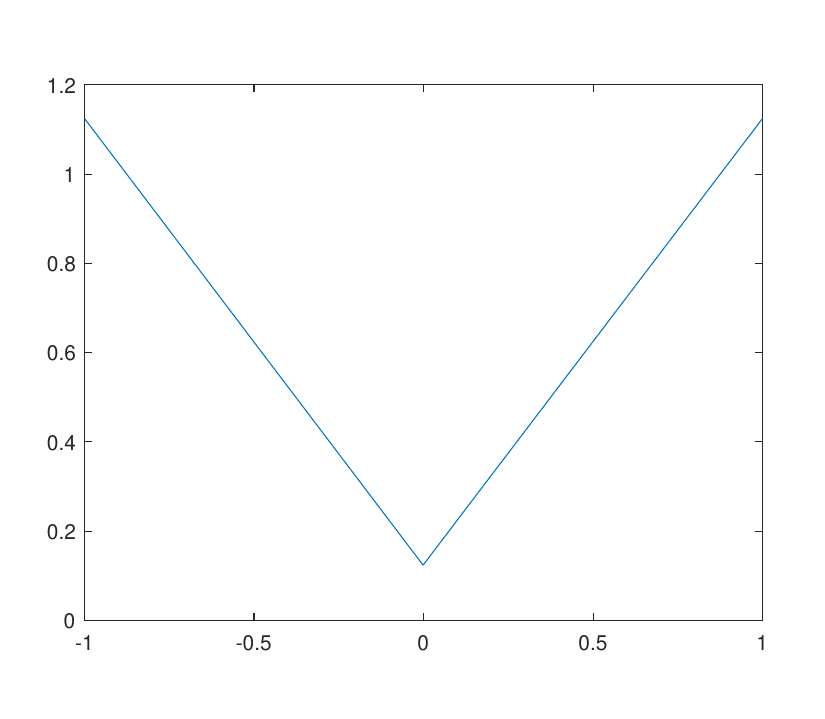}
    \caption{$f_1(t)=\sqrt{|t|}$: graph of the function, deviation and optimal spline}
    \label{fig:f1}
\end{figure}

\begin{figure}
    
    \centering
    \includegraphics[width = 0.3\textwidth]{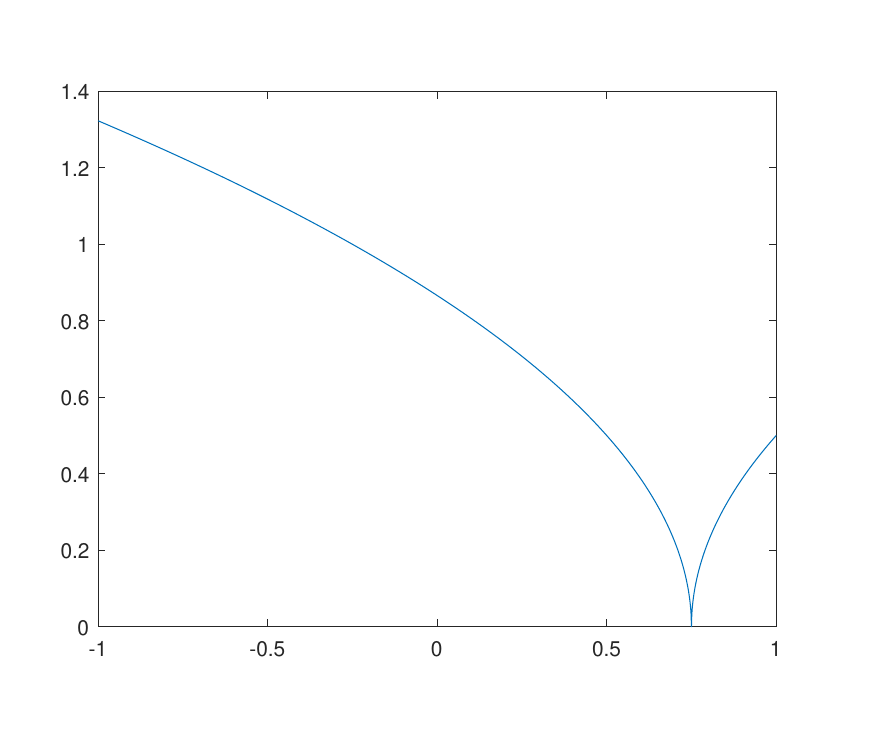}
    \includegraphics[width = 0.3\textwidth]{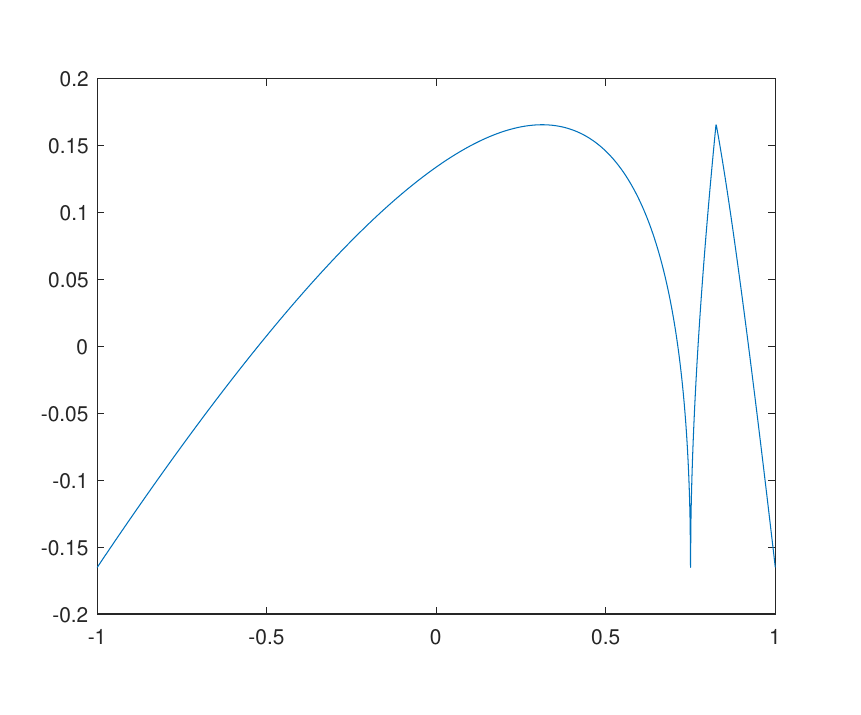}
    \includegraphics[width = 0.3\textwidth]{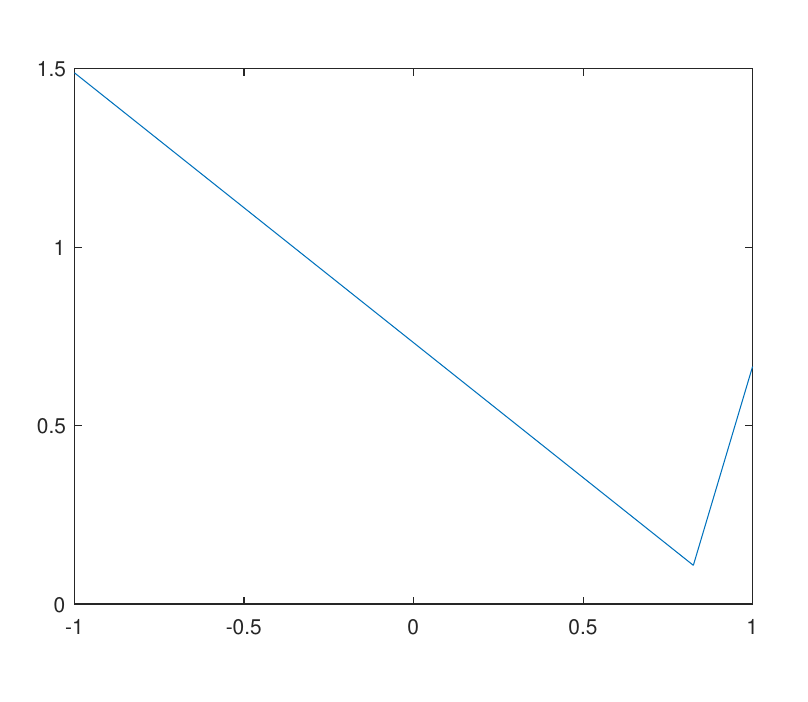}
    \caption{$f_2(t)=\sqrt{|t-0.75|}$: graph of the function, deviation and optimal spline}
    \label{fig:f2}
\end{figure}

\begin{figure}
    
    \centering
    \includegraphics[width = 0.3\textwidth]{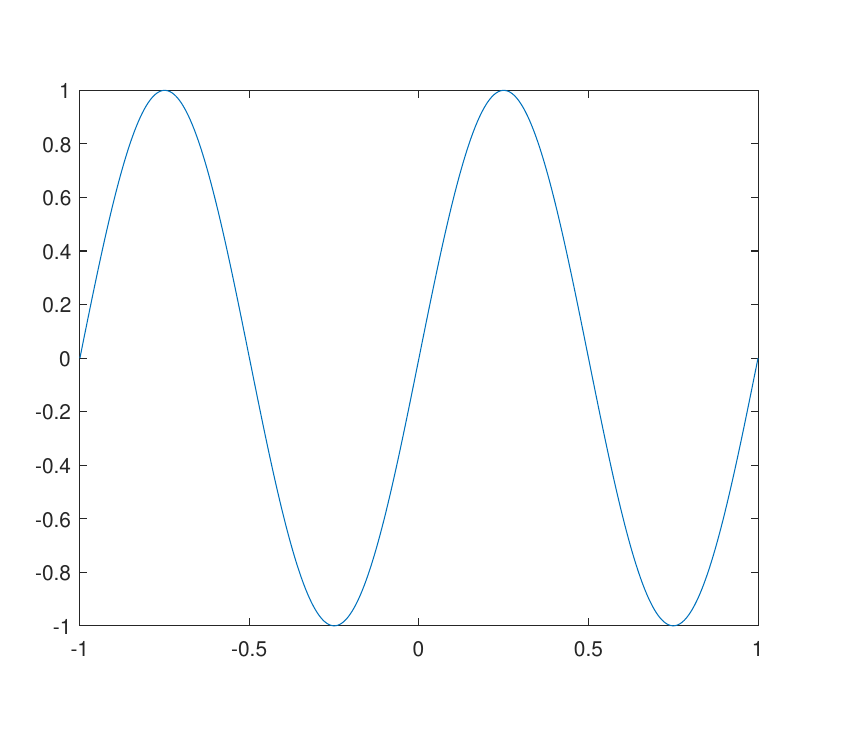}
    \includegraphics[width = 0.3\textwidth]{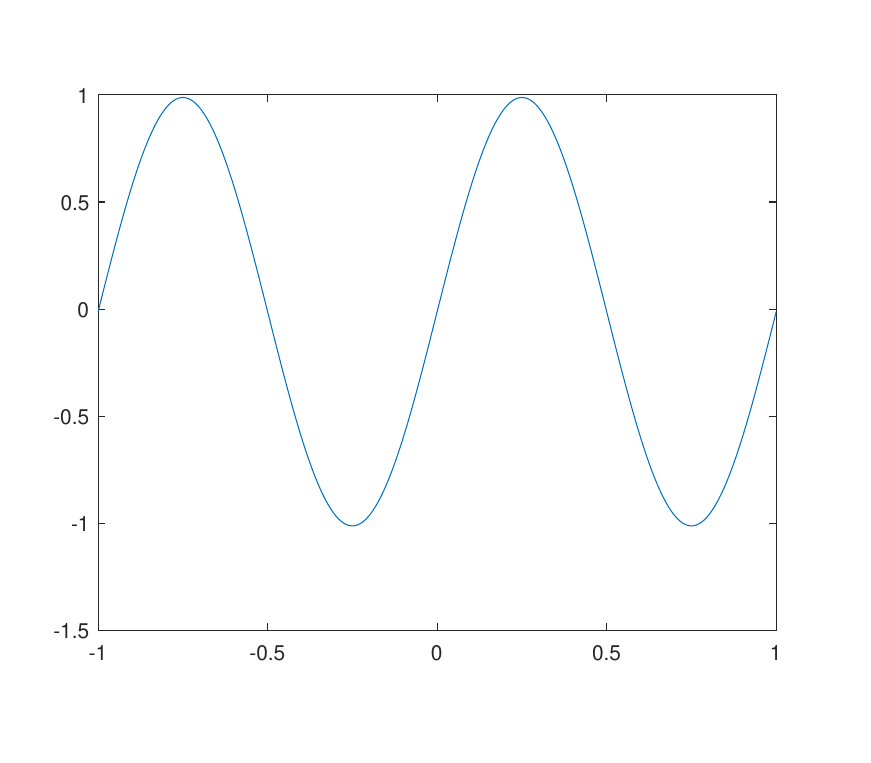}
    \includegraphics[width = 0.3\textwidth]{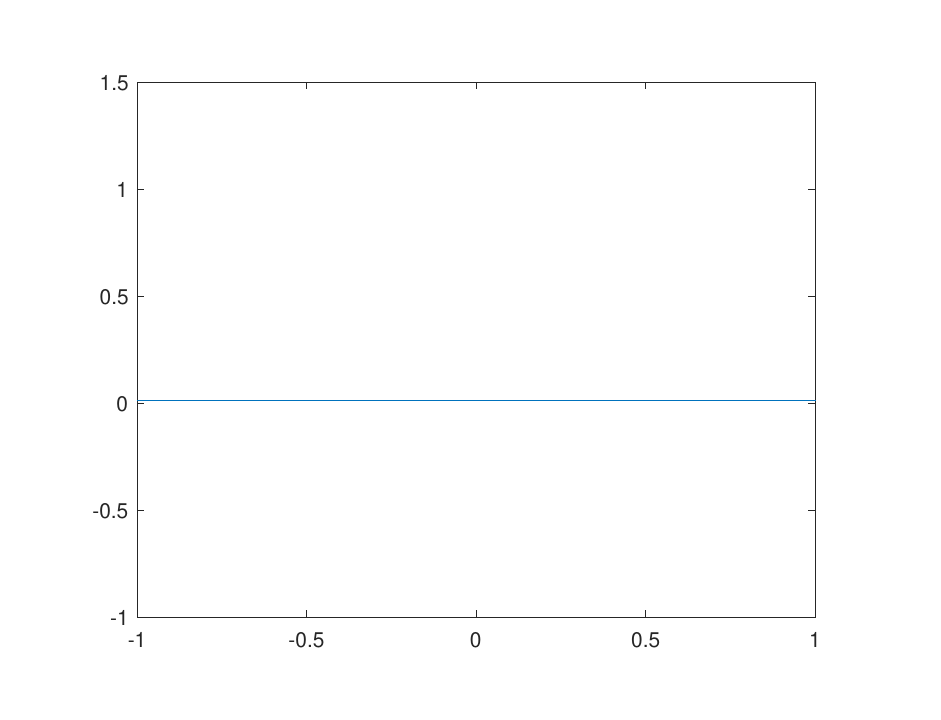}
    \caption{$f_3(t)=\sin(2\pi t)$: graph of the function, deviation and optimal spline}
    \label{fig:f3}
\end{figure}

\begin{figure}
    
    \centering
    \includegraphics[width = 0.3\textwidth]{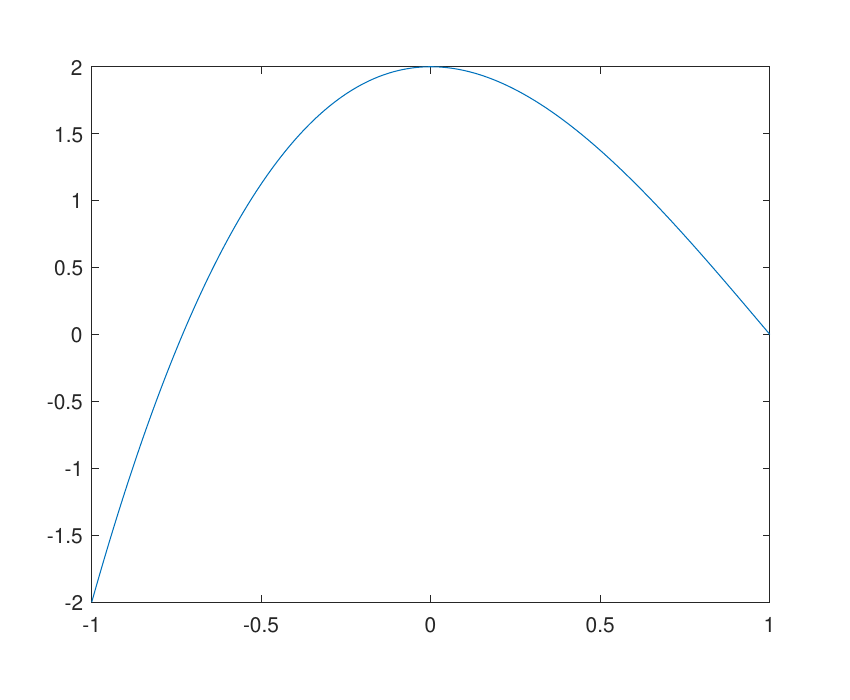}
    \includegraphics[width = 0.3\textwidth]{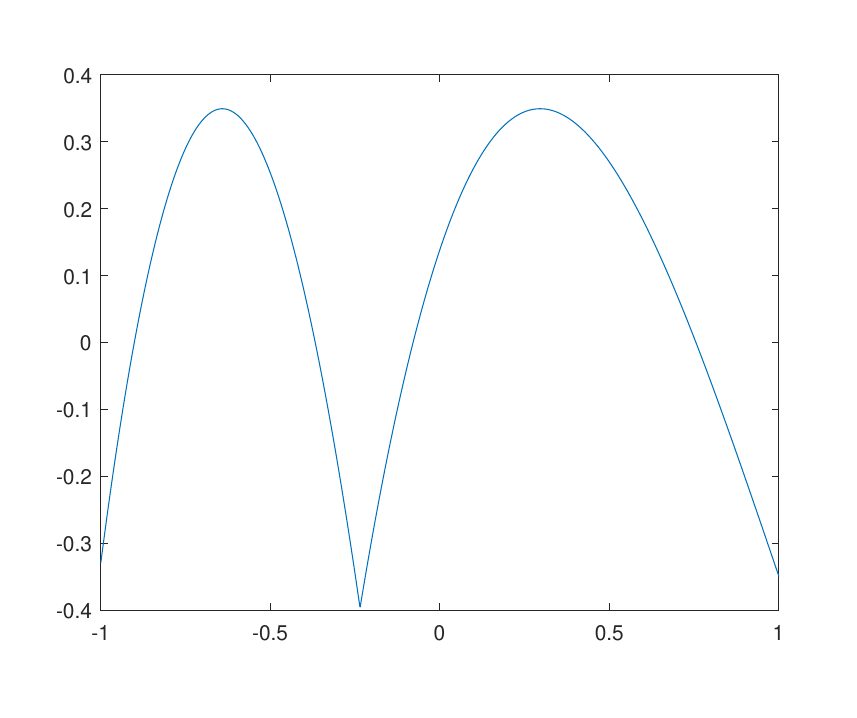}
    \includegraphics[width = 0.3\textwidth]{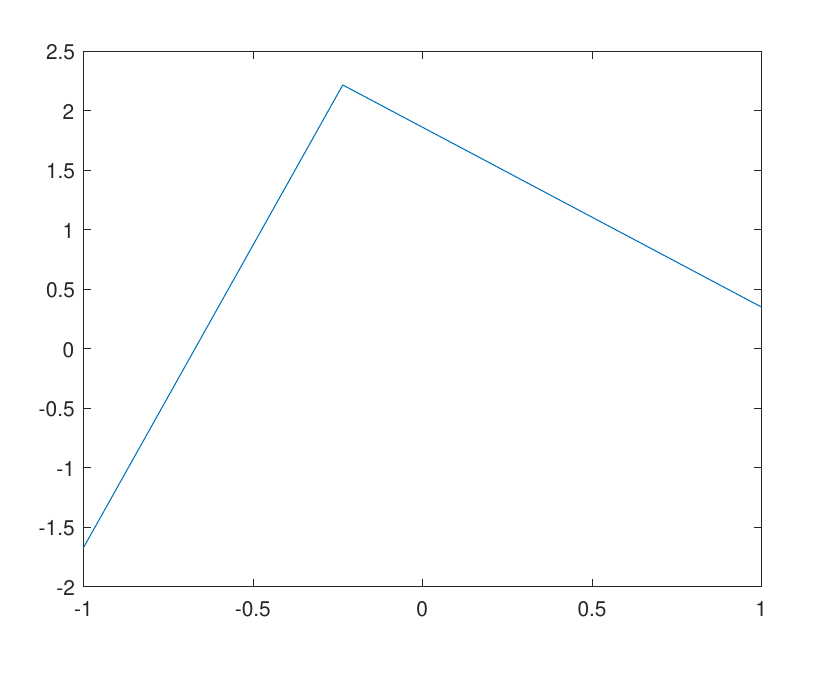}
    \caption{$f_4(t)=t^3-3t^2+2$: graph of the function, deviation and optimal spline}
    \label{fig:f4}
\end{figure}

\begin{figure}
       \centering
    \includegraphics[width = 0.3\textwidth]{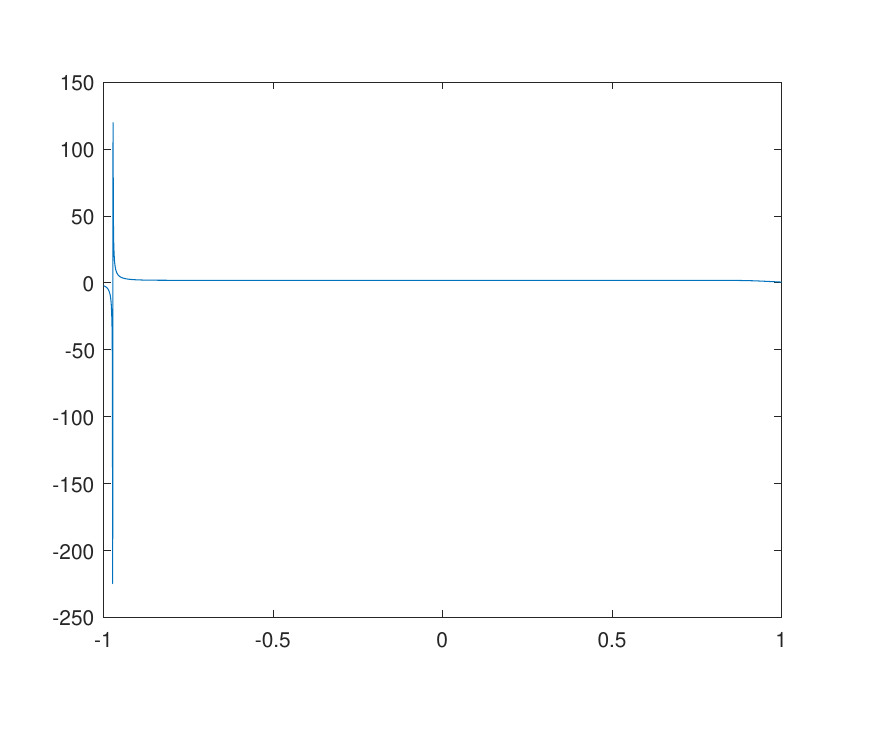}
    \includegraphics[width = 0.3\textwidth]{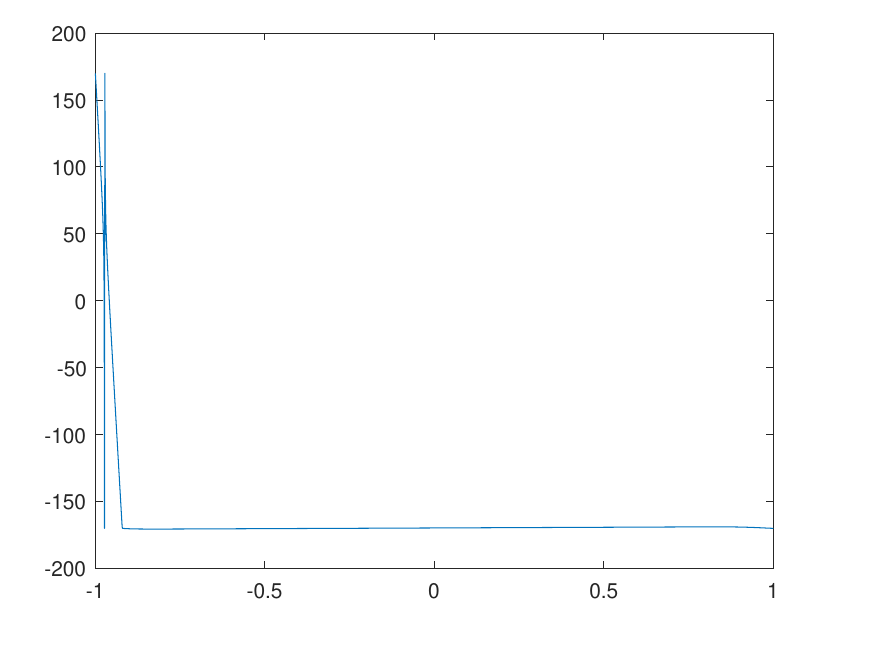}
    \includegraphics[width = 0.3\textwidth]{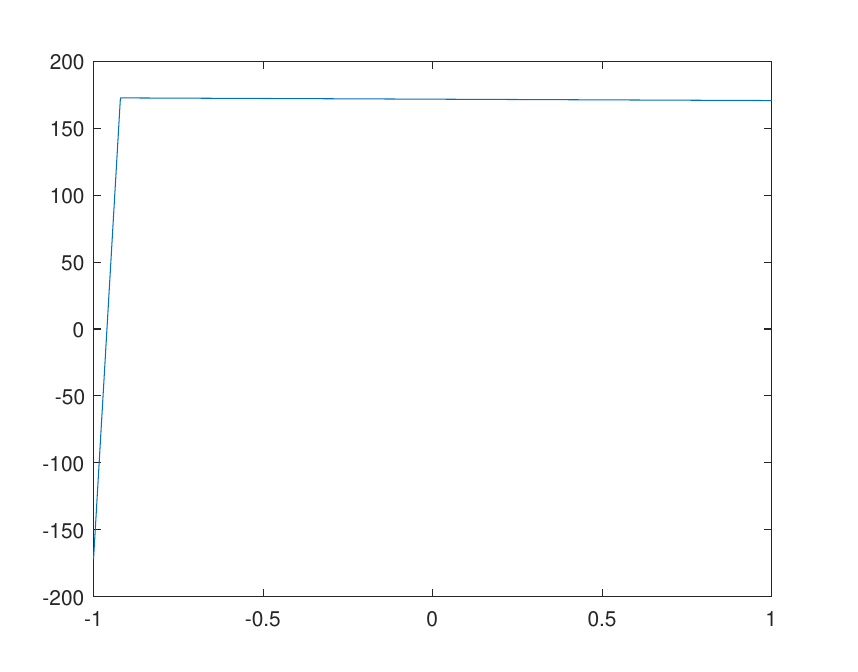}
    \caption{$f_5=1/(t^{25}+0.5)$: graph of the function, deviation and optimal spline}
     \label{fig:f5}
\end{figure}

Table~\ref{tab:res} contains the results of the numerical experiments. The first column corresponds to the function to be approximated. The second column corresponds to the parameter big-$M$ used in the corresponding mixed-linear integer programming problem. The third column is the location of the knot (N/A corresponds to the single linear piece, where the location of the knot is arbitrary and not restricted to $[c,d]$). The results for ADAM and ADAMAX (network optimisation) are in the brackets of columns~3,~4 and~6. ADAM results come first and then ADAMAX results. The fourth column is the error (maximal absolute deviation). The fifth column clarifies whether the optimal approximation corresponds to the maximum or the minimum of the linear pieces. Finally, the last column corresponds to the computational time (seconds). The results in brackets were obtained by applying neural network directly: the number of epochs for functions $f_1$, $f_2$ and $f_3$ is~50, for function $f_4$ is 100 and for function $f_5$ is 300.
\begin{table}[h]
    \centering
    \begin{tabular}{|c|c|c|c|c|c|}
    \hline
Fun  &$M$  & Knot  &Max. abs. dev. & $\max$ or $\min$ &Time (sec.)\\
\hline
      $f_1$ & 300 & 0(0.0/-0.07) &0.125(0.5/0.51) &$\max$ &332 (2.1/1.27)\\
       $f_2$ &300  & 0.75 (0.66/-0.43) &0.165(0.251/0.325) &$\max$ & 942(1.40/1.72) \\
       $f_3$ & $10^4$  &N/A(-0.48/-0.9)  & 0.999(1.044/1.002)& one piece & 35(1.23/1.32) \\
       $f_4$ & $10^4$   & -0.231(-0.17/-0.20) & 0.358(1.073/1.000)& $\min$& 365(1.41/2.06)\\
       $f_5$ & $10^5$    & -0.92(0.76/0.8) & 168.9(174/172)&$\min$ &40(2.00/3.39) \\
       \hline
    \end{tabular}
    \caption{Computational results}
    \label{tab:res}
\end{table}
For functions $f_1$, $f_2$, $f_3$ and $f_4$ the mixed-integer based models terminated at points that satisfy the sufficient optimality conditions (theorem~\ref{thm:sufficient}). In the case of the function~$f_3$, where the optimal approximation is represented by a single interval, there are four alternating points, while for the functions~$f_1$, $f_2$ and $f_4$ there are three alternating points in each subinterval (the knot $\theta$ is one of the alternating points).  In the case of the function $f_5$, the sufficient conditions are not met, but the necessary conditions from~\cite{SukhUgon2017Transactions} are met. The graphs of the functions, their linear spline approximations  and the corresponding deviation are presented in figures~\ref{fig:f1}-\ref{fig:f5}.

Comparing MILPP based models and neural networks (both, ADAM and ADAMAX), one can see that the neural network models are much faster than the optimisation models, but the maximal deviation error is lower for optimisation (more accurate results). This observation is particularly interesting in the case of function $f_4$, where the approximation of a simple smooth function appeared as a challenge for neural network-based approach, while the optimisation-based approach is accurate. If we compare ADAM and ADAMAX, the results are very similar: ADAMAX is slightly more accurate and a bit slower. The only function, where ADAM is clearly preferable is $f_2$. In our future research directions, we would like to investigate if there is something specific to neural network that makes them so fast for these particular problems.

Overall, the mixed-integer based approach for free knot linear splines with two subintervals is efficient. The choice of the parameter $M$ is crucial, but this problem is out of scope of this paper. Another interesting observation is that the neural network-based models are reasonably accurate in terms of finding the location of the knot. This can be used in the future, since one can fix the knots obtained by a neural network based model and refine the parameters, which is a convex problem and there are several efficient approaches to solve this problem.

\section{Conclusions and future research directions} \label{sec:Conclusion}

In this paper, we work with free knot linear spline approximations with at most two subintervals. This model is essential for extending to more than just one internal knot. For this model, we developed a sufficient optimality condition and an algorithm (MILPP-based approach). In our numerical experiments, we use five functions. All these functions are challenging for linear spline approximation with one internal knot: abrupt changes, nonsmoothess, periodic nature, non-Lipschitzness.  In all the experiments, except the last one (function~$f_5$) the approximations satisfy the sufficient optimality condition developed in this paper. Function~$f_5=1/(t^{25}+0.5)$ is a very complex function to approximate, but the final approximation satisfies the strongest known necessary optimality condition from~\cite{SukhUgon2017Transactions}.

It is also important to note that the proposed MILPP-based approach is also suitable for approximating multivariate functions. In this case, the size of the discretisation lattice (total number of discretisation points) is increasing and therefore the dimension of the corresponding optimisation problems (due to $c_i$) is increasing. 

Our future research directions include the following.
\begin{enumerate}
    \item The development of a Remez-type algorithm for linear spline approximation, where the point update procedure relies on the MILPP model proposed in this paper.
    \item The investigation how the results obtained by a fast neural network model can be used as an initial guess for free knot linear spline approximation. Moreover, it is important to understand what is so specific in modern neural network models that make them reasonably efficient at finding the location of the knots and if it will still be the case for more than one internal knot. 
\end{enumerate}
\section*{Acknowledgement}
This research was supported by the Australian Research Council (ARC),  Solving hard Chebyshev approximation problems through nonsmooth analysis (Discovery Project DP180100602).


\begin{thebibliography}{99}

\bibitem{Nurnberger}
G. Nürnberger, Approximation by spline functions, Springer-Verlag, 1989.

\bibitem{Schum}
L. Schumaker, Uniform approximation by {C}hebyshev spline functions. II: free knots, SIAM Journal of Numerical Analysis 5 (1968), 647–656.

\bibitem{OpenProblems} P. Borwein, I. Daubechies, V. Totik, and G. Nürnberger, Bivariate segment approximation
and free knot splines: Research problems 96-4, Constructive Approximation 12 (1996), no. 4, 555–558

\bibitem{SukhUgon2017Transactions} Nadezda Sukhorukova and Julien Ugon, Characterisation theorem for best polynomial spline approximation with free knots, 
Trans. Amer. Math. Soc. 369 (2017), 6389-6405

\bibitem{CrouSukhUgon} J. P. Crouzeix, N. Sukhorukova, and J. Ugon, Characterization theorem for best polynomial spline
approximation with free knots, variable degree and fixed tails, J. of Optimization Theory and Applications, Vol. 172, No. 3, 950–964 (2017)

\bibitem{NurnSchuSomm} G. Nürnberger, L. Schumaker, M. Sommer, and H. Strauss, Approximation by generalized splines, Journal of Mathematical Analysis and Applications 108 (1985), 466–494 

\bibitem{Cheb} Chebyshev, P., Théorie des mécanismes connus sous le nom de parallélogrammes,   {Mémoires des Savants étrangers présentés à l’Académie de Saint-Pétersbourg}, vol.  {7} (1854),  {539--586}. 

\bibitem{Chen} Chen K. K. (2016). The Upper Bound on Knots in Neural Networks. Working paper arXiv:1611.09448.

\bibitem{Goodfellow2016} Goodfellow, I., Bengio, Y., Courville, A.: Deep Learning. MIT Press (2016). http://www.deeplearningbook.org

\bibitem{Sun2020OptimDeepLearning} Sun, R.Y.: Optimization for deep learning: An overview. Journal of the Operations Research Society of China 8, 249–294 (2020)

\bibitem{Vidal2017} Haeffele, B.D., Vidal, R.: Global optimality in neural network training. In: 2017 IEEE Conference on Computer Vision and
Pattern Recognition, CVPR 2017, Honolulu, HI, USA, July 21-26, 2017, pp. 4390–4398. IEEE Computer Society (2017). doi
10.1109/CVPR.2017.467. https://doi.org/10.1109/CVPR.2017.46

\bibitem{Cybenko} Cybenko, G.: Approximation by superpositions of a sigmoidal function. Mathematics of Control, Signals and Systems 2,
303–314 (1989

\bibitem{Hornik1991} Hornik, K.: Approximation capabilities of multilayer feedforward networks. Neural networks 4(2), 251–257 (1991)

\bibitem{Pinkus1993} Leshno, M., Lin, V.Y., Pinkus, A., Schocken, S.: Multilayer feedforward networks with a nonpolynomial activation function
can approximate any function. Neural Networks 6(6), 861–867 (1993). https://doi.org/10.1016/S0893-6080(05)80131-

\bibitem{Pinkus1999} Pinkus, A.: Approximation theory of the MLP model in neural networks. Acta Numerica 8, 143–195 (1999). doi
10.1017/S0962492900002919

\bibitem{Kolmogorov1957} Kolmogorov, A.N.: On the representation of continuous functions of many variables by superposition of continuous functions
of one variable and addition. Dokl. Akad. Nauk SSSR 114, 953–956 (1957)

\bibitem{Arnold1963} Arnold, V.: On functions of three variables. Dokl. Akad. Nauk SSSR 114, 679–681 (1957). English translation: Amer. Math.
Soc. Transl., 28 (1963), pp. 51–54

\bibitem{PeirisRoshchinaSukhorukova2023}  Peiris, V., Roshchina, V., Sukhorukova, N.: Artificial Neural Networks with uniform norm based loss function. In PRESS, accepted by Advances in Computational Mathematics in 2023. 

\bibitem{hansson2017} Hansson, M., and Olsson, C.: Feedforward neural networks with ReLU activation functions are linear splines. Bachelor's Theses in Mathematical Sciences (2017).

\end{thebibliography}
\end{document}